\documentclass[journal]{IEEEtran}
\usepackage[T1]{fontenc}
\usepackage{times}
\usepackage{amsmath,amsfonts,amssymb}
\usepackage{amsthm}
\usepackage{bbm}
\usepackage{algorithmic}
\usepackage{algorithm}
\usepackage{array}
\usepackage[caption=false,font=normalsize,labelfont=sf,textfont=sf]{subfig}
\usepackage{textcomp}
\usepackage{stfloats}
\usepackage{url}
\usepackage{verbatim}
\usepackage{graphicx}
\usepackage{cite}
\usepackage{hyperref}
\usepackage[table]{xcolor}
\usepackage{float}
\usepackage{tabularx}
\usepackage{multirow}
\usepackage{bbold}
\usepackage{dsfont}
\usepackage{enumitem}
\usepackage{newtxtext}
\usepackage{pgfplots}
\pgfplotsset{compat=1.18}
\usepackage{booktabs, siunitx}
\newcolumntype{Y}{>{\centering\arraybackslash}X} 

\DeclareMathOperator{\hyp}{hyp}
\DeclareMathOperator{\conv}{conv}

\begin{document}

\title{GPU-Accelerated Dynamic Programming for Multistage Stochastic Energy Storage Arbitrage}

\author{Thomas Lee, \IEEEmembership{Student Member, IEEE}, Andy Sun, \IEEEmembership{Senior Member, IEEE.}

\thanks{Thomas Lee is with the Institute for Data, Systems, and Society, MIT.}
\thanks{Andy Sun is with the Sloan School of Management, MIT.}}

\markboth{}%
{}



\maketitle

\begin{abstract}
We develop a GPU-accelerated dynamic programming (DP) method for valuing, operating, and bidding energy storage under multistage stochastic electricity prices. Motivated by computational limitations in existing models, we formulate DP backward induction entirely in tensor-based algebraic operations that map naturally onto massively parallel GPU hardware. Our method accommodates general, potentially non-concave payoff structures, by combining a discretized DP formulation with a convexification procedure that produces market-feasible, monotonic price-quantity bid curves. Numerical experiments using ISO-NE real-time prices demonstrate up to a 100x speedup by the proposed GPU-based DP method relative to CPU computation, and an 8,000x speedup compared to a commercial MILP solver, while retaining sub-0.3\% optimality gaps compared to exact benchmarks. 

\end{abstract}


\section{Introduction}

Energy storage plays a growing role in modern power systems. Energy arbitrage, which shifts energy between lower and higher-valued periods, is a critical revenue source for storage projects. Increased grid volatility, driven by load growth and renewables integration, raises computational challenges for accurate energy storage modeling.

Since dispatch decisions must be made before prices are known, deterministic perfect-foresight models, such as typical LP formulations, can significantly overestimate arbitrage value \cite{sioshansi2021energy}. Accurately modeling non-simultaneous charging and discharging remains challenging: LP relaxations cannot fully eliminate complementarity violations since the convex hull may involve exponentially many cuts \cite{pozo2023convex}, and the issue may become more severe as negative prices increase in frequency due to renewables deployment and grid congestion.

Further, electricity markets rely on \emph{price-quantity} bids instead of solely self-scheduled quantities. Prior works tackling the storage bid formation problem use approximate dynamic programming \cite{jiang2015optimal} or binary variables to indicate bid-clearing status \cite{krishnamurthy2017energy}, but these models only allow one bid segment per time period rather than more general bid curves. A recent analytic dynamic programming (DP) approach constructs energy arbitrage bid curves from value-function subgradients and reports significant speedups compared to stochastic dual dynamic programming (SDDP) \cite{xu2020operational}, but it requires linear payoffs, prohibits discharging at negative prices, and still relies on state space discretization for numerical implementation (concave payoffs require discrete actions as well).

With theoretical roots in DP, deep reinforcement learning (RL) offers fast inference due to GPU-accelerated neural networks and has been applied to stochastic energy arbitrage \cite{cao2020deep}, but it requires substantial upfront training effort, and does not generate market-compatible bid curves.

In contrast to using GPU computation for neural networks (e.g. for RL or ML-based optimization surrogates), recent research efforts have applied GPU acceleration directly within novel optimization solvers, including for linear programming (LP) \cite{lu2025cupdlp} and nonlinear programming (NLP) \cite{shin2024accelerating}. Inspired by these recent advances that utilize the optimization algorithm structure to leverage GPU hardware, in this work we develop a discretized DP method with a matrix-tensor structure that naturally enables efficient GPU calculations.

\subsection{Contributions}
This work makes the following contributions:
\begin{itemize}
    \item We design a tensor-based backward induction algorithm to solve a dynamic programming formulation of the multistage stochastic energy arbitrage problem. 

    \item We implement our algorithm to leverage GPU hardware, achieving up to 100x speedup versus CPU baselines.

    \item We propose a general convexification approach to derive storage bidding curves under general payoff functions.
    
    \item Using a test case with extended negative prices, we demonstrate the advantage of the discretized DP method over existing approximate approaches, and show an 8,000x speedup over a commercial MILP solver. 
        
    \item We use our DP method to quantify the economic value of stochastic modeling and price-quantity bidding.
    
\end{itemize}

\section{Models}
\label{sec:models}
\subsection{Deterministic energy arbitrage}
The deterministic, price-taker energy arbitrage problem is
\begin{subequations}
\label{eq:arbitrage}
    \begin{align}
        \max\nolimits_{\{p_t,s_t\}_{t=1}^T} \ & \sum\nolimits_{t=1}^T \lambda_t p_t
        \\
        \text{s.t.} \ & 
        -\overline{p} \leq p_t \leq \overline{p}, \quad \forall t\in [T], 
        \label{eq:feas:power}
        \\
         & 0 \leq s_t \leq \overline{s}, \quad \forall t \in [T],
         \label{eq:feas:energy}
         \\
         & s_{t} = s_{t-1} + F(p_t), \quad \forall t \in [T], \label{eq:feas:soc}
    \end{align}
\end{subequations}
which, given an initial energy level $s_0$ along with power capacity $\overline{p}$ and energy capacity $\overline{s}$, decides a sequence of net power outputs $\{p_t\}$ and energy levels $\{s_t\}$, to maximize the total storage profits which are linear to the power prices $\{\lambda_t\}$. The \emph{state transition} function in Eq. \eqref{eq:feas:soc} is
\begin{align}
    F(p_t) := -\begin{cases}
        \frac{1}{\eta}p_t, &\text{if } p_t \geq 0, \text{ (discharge)}
        \\
        \eta p_t, &\text{if } p_t < 0 \ \text{ (charge)}
    \end{cases} \label{eq:F}
\end{align}
which describes the power-to-energy conversion with an efficiency factor $\eta \in (0,1]$, enforcing charge-discharge complementarity for battery energy storage systems. 

Prior works typically formulate Eq. \eqref{eq:arbitrage}-\eqref{eq:F} with separate charging $p^c_t$ and discharging $p^d_t$. An exact MILP reformulation of Eq. \eqref{eq:arbitrage}-\eqref{eq:F} is to combine Eq. \eqref{eq:arbitrage} with the added constraints
\begin{subequations}
\label{eq:MILP}
\begin{align}
    p_t &= p^d_t - p^c_t,\quad \forall t\in [T] \label{eq:MILP:3a}
    \\
    F(p_t) &= p^c_t \eta - p^d_t/\eta,\quad \forall t\in [T] \label{eq:MILP:3b}
    \\
    0 &\leq p^c_t \leq \overline{p} z_t, \quad \forall t \in [T], 
    \\
    0 &\leq p^d_t \leq \overline{p} (1-z_t),\quad \forall t\in [T] \label{eq:MILP:3c}
    \\
    z_t &\in \{0,1\}, \quad \forall t \in [T]. \label{eq:MILP:3d}
\end{align}
\end{subequations}
This is standard in the literature \cite{pozo2023convex}. Note $[T] \equiv \{1,...,T\}$.

\subsection{Multistage stochastic energy arbitrage}
In the stochastic setting, we model power prices as \emph{stagewise independent} random variables, following \cite{xu2020operational}. That is, the random prices $\{\lambda_t\}_t$ are mutually independent, each with a probability distribution. Given a feasible incoming state $s_{t-1} \in [0, \overline{s}]$, the set of feasible power actions is the continuous interval
\begin{align}
    \mathcal{P}(s_{t-1}) = \left[
    {-}\min\{\overline{p}, (\overline{s} - s_{t-1})/\eta \}, \ \min\{ \overline{p}, s_{t-1} \eta \}
    \right], \label{eq:feas_power}
\end{align}
because when $p_t\geq 0$, we need $p_t \leq \overline{p}$ from \eqref{eq:feas:power} and $0 \leq s_{t-1} - \frac{1}{\eta}p_t$ from \eqref{eq:feas:soc}-\eqref{eq:F}. Conversely, when $p_t < 0$, we need $-\overline{p} \leq p_t$ and $s_{t-1} - \eta p_t \leq \overline{s}$. Thus Eq. \eqref{eq:feas_power} exactly encapsulates the constraints \eqref{eq:feas:power}-\eqref{eq:F} for two adjacent periods $t{-}1$ and $t$. 

Assuming stagewise independence, the dynamic programming problem can be written, adapting \cite{xu2020operational}'s formulation, using the \emph{value function} defined recursively based on Bellman optimality for each $t \in [T]$:
\begin{align}
    Q_{t-1}(s_{t-1}, \xi_t) = 
    \max\limits_{p_t \in \mathcal{P}(s_{t{-}1}) } \{ 
         \lambda_t p_t + V_{t}\left(s_{t{-}1} + F(p_t)\right)\},
    \label{eq:storage_DP}
\end{align}
where the uncertainty here is just in prices $(\xi_t = \lambda_t)$. Using the optimality condition for Eq. \eqref{eq:storage_DP}, we will be able to find the price-quantity bid curve in Section \ref{sec:bidding}.

The \emph{expected} value function is defined recursively:
\begin{align}
        V_{t-1}(s_{t-1}) := \mathbbm{E}_{\xi_t} \left[Q_{t-1}(s_{t-1}, \xi_{t})\right], \quad \forall t \in \{0\} \cup [T], \label{eq:V}
\end{align}
with the final-period base case, describing the value of any residual state-of-charge, defined as $V_{T}(s_{T}) := 0, \ \forall s_T \in \mathbbm{R}.$ Altogether, Eq. \eqref{eq:storage_DP}-\eqref{eq:V} constitute the dynamic programming (DP) problem for multistage stochastic energy arbitrage. The objective expression of Eq. \eqref{eq:storage_DP} embeds the state transition of \eqref{eq:feas:soc}. The multistage problem \eqref{eq:storage_DP} solves for the optimal policies of storage power dispatch $p_t$, which are contingent upon the incoming energy level $s_{t-1}$ as well as the price realization $\lambda_t$.

\subsection{Price-quantity function bidding}
\label{sec:bidding}
In order to participate in real-time (RT) electricity markets, resources including energy storage submit \emph{price-quantity function bids}, consisting of \emph{price-quantity pairs} that enter the grid operator's economic dispatch problem. We assume the bid curve for delivery during time $t$ should be calculated prior to time $t$, following \cite{xu2020operational}. Define the composed function
\begin{align}
    U_t(p_t; s_{t-1}) := V_t(s_{t-1} + F(p_t)),
\end{align}
which maps action $p_t$ to the value obtained at the next stage. 

Consider the setting where $V_t(s_t)$ is concave and monotonic non-decreasing in its energy level argument $s_t$. The concave $V_t(\cdot)$ is composed with $F(p_t)$, which preserves the concavity of $U(p_t;s_{t-1})$ in its argument $p_t$ (see p.32 in \cite{rockafellar1970convex}).
Then, given any $s_{t-1}$, the optimality condition for Eq. \eqref{eq:storage_DP} is
$\lambda_t 
    \in 
      -\partial U_t(p_t; s_{t-1}), \label{eq:opt1}
$
which depends on the superdifferential of the concave $U_t(p_t;s_{t-1})$ with respect to $p_t$. Then
\begin{align}
    b(p_t;s_{t-1}) := -\partial {U}_t(p_t; s_{t-1}), \label{eq:bids}
\end{align}
is a set-valued function that maps a feasible power output $p_t \in \mathcal{P}(s_{t-1})$ to a corresponding set of marginal costs, within which $p_t$ is optimal. The concavity of $U_t(\cdot;s_{t-1})$ means $b(p_t)$ is monotonic non-decreasing in $p_t$.
The overall optimal bidding curve is formed by the graph
\begin{align}
(
  \mathcal{P}(s_{t-1}),
  \ 
  b(\mathcal{P}(s_{t-1});s_{t-1}) 
). \label{eq:bid_curves}
\end{align}

Critically, the bid curve written in Eq. \eqref{eq:bid_curves} must be decided during interval $t-1$ and corresponds to physical delivery during interval $t$. From Eq. \eqref{eq:bids}, this bid curve only depends on information embedded in $V_t(\cdot)$, which is an expectation over the $\lambda_t$ and remains uncertain as of time $t-1$. Thus the bid curve construction is non-anticipative and is feasible to the real-life imperfect information structure.

In general, $U_t(\cdot,s_{t-1})$ could be non-concave. For example, if prices are sufficiently negative then $V_t(s_t)$ could actually decrease as $s_t$ increases, i.e. when having more headroom is beneficial because charging is profitable. In this case, the function could be \emph{convexified} first, to $\tilde{U}_t(p_t;s_{t-1})$, such that
\begin{align}
    \hyp \tilde{U}_t(p_t;s_{t-1}) = \conv (\hyp {U}_t(p_t;s_{t-1}) ),
\end{align}
where $\hyp$ denotes the hypograph of a function and $\conv$ denotes the convex hull. This step can be performed numerically, for example with the Graham scan \cite{Graham1972ConvexHull}. Then, Eq. \eqref{eq:bids} can be modified to produce a \emph{monotonized} bid curve using
\begin{align}
    \tilde{b}(p_t;s_{t-1}) := - \partial \tilde{U}_t(p_t; s_{t-1}),
    \label{eq:bids_convex}
\end{align}
which now meets the monotonicity required by market rules. The next section describes the numerical implementation of the overall dynamic programming algorithm.

\section{Dynamic programming algorithm}

\subsection{Discretization}
\label{sec:discretization}
Assume that the state space is discretized to uniform multiples of a divisible scalar $\delta > 0$, such that $n^s = \overline{s}/\delta \in \mathbbm{N}$. Explicitly write out the vector of these state levels as
\begin{align}
    \hat{s} := \left[0, \delta, 2\delta, ..., \overline{s}  \right]^\top \in \mathbbm{R}^S,
    \quad \quad
    S := n^s + 1.
\end{align}

Define the discretized power actions as the vector
\begin{subequations}
\label{eq:discrete_actions}
\begin{align}
    \hat{p} :=& \left[
        -\overline{p}
    ,
    ...,
    -2{\delta}/{\eta}, 
    -{\delta}/{\eta}, \ 0, \ 
    \delta\eta,
    2\delta\eta,
    ...,
    \overline{p}
    \right]^\top,
    \\
    P :=& \dim(\hat{p}) = n^c + n^d + 1,
    \\
    n^c :=& \left\lceil \overline{p}\eta  / \delta  \right\rceil, \quad n^d = \left\lceil \overline{p} / (\delta \eta) \right\rceil,
\end{align}
\end{subequations}
where $n^c, n^d$ are the number of charge and discharge actions, plus the zero-power action. The two $\{-\overline{p},\overline{p}\}$ endpoints are deliberately included based on the empirical observation that optimal storage dispatch very often occurs at the $\{-\overline{p}_t, 0, \overline{p}_t\}$ levels (which is reinforced by the theoretical existence of extreme point LP solutions).

Finally, assume that each random variable $\lambda_t$ is modeled with a discrete probability mass function with a finite support (e.g. by randomly sampling from a continuous distribution of $\lambda_t$), accessible as a vector $\hat{\lambda}_t$ of price levels along with a vector $\hat{\pi}_t$ of their discrete probabilities,
\begin{subequations}
\begin{align}
    \hat{\lambda}_t,\ \hat{\pi}_t &\in \mathbbm{R}^R, \quad \forall t\in[T],
    \\
    \hat{\pi}_{tr} &= \mathbbm{P}(\lambda_t = \hat{\lambda}_{tr}),
    \quad \forall t\in[T], r \in [R]
\end{align}
\end{subequations}
where $R$ is the cardinality of the finite discrete support of $\lambda_t$, i.e. the number of distinct price levels used to model its distribution, and such that $\sum_{r=1}^R \hat{\pi}_{tr} =1$ for all $t \in [T]$.

\begin{figure*}[hb]
    \centering
\includegraphics[width=0.99\textwidth]{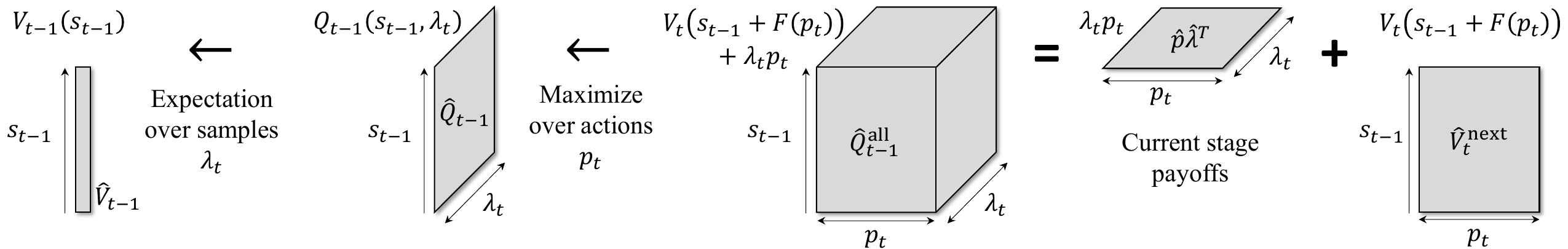}
    \caption{Diagram illustrating Algorithm \ref{alg:algorithm}. Vector / matrix / tensor objects from Alg. \ref{alg:algorithm} are illustrated as gray shapes. Labels at the top correspond to function names from Section \ref{sec:models}, and each axis is labeled with its corresponding variable name. The ``+'' operation here represents the tensor sum defined in Line \ref{eq:Qhat_all}.
    }
    \label{fig:diagram}
\end{figure*}

\subsection{Backward induction as matrix algebra}

We now impose the discretization scheme of Section \ref{sec:discretization} onto the DP problem of Eq. \eqref{eq:feas_power}-\eqref{eq:V}. The resulting discretized problem is thus a restriction of the original DP, since all actions and states are still feasible. Algorithm \ref{alg:algorithm} solves the discretized DP using purely matrix-algebraic steps. 

\begin{algorithm}[ht]
\caption{Tensor-based DP backward induction}
\label{alg:algorithm}
\begin{algorithmic}[1]
    \REQUIRE $\hat{s} \in \mathbbm{R}^S, \ \hat{p} \in \mathbbm{R}^P, \ 
    \{\hat{\lambda}_t, \hat{\pi}_t\}_{t=1}^T$, where $\hat{\lambda}_t,\hat{\pi}_t \in \mathbbm{R}^R$.

    \ENSURE Expected value vectors $\{\hat{V}_t\}_{t=0}^T$.
    
    \STATE $\hat{V}_T \gets \mathbf{0}_S$. \hfill (\emph{base case})

    \STATE ${\sigma} \gets \hat{s} (\mathbf{1}_P)^\top + (\mathbf{1}_S) F(\hat{p})^\top$.
    \label{eq:next_state_levels}
    \hfill (\emph{next-state levels})
    
    \STATE $z \gets \frac{1}{\delta} \min\{0, \max\{{\sigma} , \overline{s} \} \} + 1$. \label{eq:next_state_indices}
    \hfill (\emph{next-state proxies})

    \STATE $z^- \gets \lfloor z \rfloor, \quad
    z^+ \gets \lceil z \rceil$.    \label{eq:next_state_indices_discrete}
    \hfill
    (\emph{next-state indices})

    \STATE \label{eq:interpolation_factors}
    $b \gets ({z} - {z}^-) /  ({z}^+ - {z}^-)$. 
    \hfill (\emph{interpolation factors})

    \FOR{$t = T, ..., 0$}
    
     \STATE $\hat{V}_t^{\mathrm{next}} \gets (1 - b)\odot \hat{V}_{t}[z^-] + b\odot \hat{V}_{t}[z^+]$. 
    \label{eq:Vnext}

     \STATE $\hat{V}^\mathrm{next}_t [\mathbbm{1}(\sigma < 0) + \mathbbm{1}(\sigma > S)] \gets -\infty$. \label{eq:infeasibility}\hfill (\emph{infeasibility}) 

     \STATE $    \hat{Q}^\mathrm{all}_{t{-}1} \gets \mathbf{1}_S \otimes \hat{p} \hat{\lambda}_t^\top + \hat{V}_t^{\mathrm{next}} \otimes \mathbf{1}_R^\top$. \label{eq:Qhat_all}

    \STATE $\hat{Q}_{t{-}1}[i,r] \gets \max\limits_{j\in[P]} \hat{Q}^\mathrm{all}_{t-1}[i,j,r], \quad \forall i\in[S],r\in[R]$. \label{eq:Qhat}

    \STATE $\hat{V}_{t-1} \gets \hat{Q}_{t-1} \hat{\pi}_{t}$. \hfill (\emph{expectation}) \label{eq:expectation}

    \ENDFOR
\end{algorithmic}
\end{algorithm}

A key part of the algorithm design is efficiently evaluating the next-stage value functions. This is achieved by Lines \ref{eq:next_state_levels}-\ref{eq:interpolation_factors}, which pre-compute relevant vector indices and interpolation factors as constants across all $t$. Line \ref{eq:next_state_levels} performs the ``outer sum'' on the state and power vectors to find $\sigma\in \mathbbm{R}^{S\times P}$, consisting of all pairwise sums from $\hat{s}$ and $\hat{p}$. Line \ref{eq:next_state_indices} projects $\sigma$, elementwise, onto $[0,\overline{s}]$, then transforms it to $z \in \mathbbm{R}^{S\times P}$ with elements in the continuous range $[1,S]$, forming continuous ``proxies'' of the next-state indices. Line \ref{eq:next_state_indices_discrete} converts $z$ into matrices of integer indices $z^-,z^+ \in [S]^{S \times P}$, where floor and ceiling operations are applied elementwise. Then Line \ref{eq:interpolation_factors} calculates interpolation factors (with slight abuse of notation, where division applies elementwise). Note that
$$F(\hat{p}) = [\overline{p}\eta, ..., 2\delta, \delta, 0, -\delta, -2\delta, ..., \overline{p}/\eta]^\top,$$
which guarantees that all of $\sigma$'s columns (except the first and last) consist of combining integer multiples of $\delta$. This means columns 2 to $P{-}1$ in matrix $z$ all contain integer elements in set $[S]$, associated with 0 interpolation factor within $b$. This choice of ``recombining'' states and power actions helps to reduce any numerical interpolation errors, since we would only need to interpolate for the extreme power endpoints $\{-\overline{p}, \overline{p}\}$. 

Next, Algorithm \ref{alg:algorithm} performs backward induction. At each stage $t$, Line \ref{eq:Vnext} uses the index matrices $\hat{z}^-,\hat{z}^+$ to extract the appropriate elements of $\hat{V}_t$ and applies linear interpolation for the first and last columns (corresponding to $-p_t$ and $p_t$). In Line \ref{eq:Vnext}, brackets denote indexing, and $\odot$ denotes the elementwise Hadamard product. This produces a tentative next-state value matrix $\hat{V}^\mathrm{next}_t \in \mathbbm{R}^{S\times P}$, which is a discretization of $V_t(s_{t-1}+F(p_t))$. Line \ref{eq:infeasibility} then assigns $-\infty$ to the infeasible (state, power) combinations, in $\hat{V}^\mathrm{next}_t$, to ensure they are not chosen during maximization. The outer product matrix $\hat{p} \hat{\lambda}_t^\top \in \mathbbm{R}^{P\times R}$ describes the stage-specific payoffs, and the tensor $\hat{Q}^\mathrm{all}_{t-1} \in
\mathbbm{R}^{S\times P \times R}$ is created in Line \ref{eq:Qhat_all} by broadcasting the matrices with tensor products $\otimes$ onto common dimensions.

At this point, $\hat{Q}^\mathrm{all}_{t-1}$ corresponds to the inner expression within Eq. \eqref{eq:storage_DP}'s maximization objective. Line \ref{eq:Qhat} then maximizes over the actions, resulting in the matrix $\hat{Q}_{t-1} \in \mathbbm{R}^{S\times R}$. Next, Line \ref{eq:expectation} computes the expectation over the discrete random variable $\lambda_t$; this is done as a matrix-vector product, producing $\hat{V}_{t-1} \in \mathbbm{R}^S$. Then the algorithm inductively proceeds backward, as illustrated in Fig. \ref{fig:diagram}. Finally, the bidding function in Eq. \eqref{eq:bids} or \eqref{eq:bids_convex} can be evaluated using numerical finite differences based on the $\hat{V}_t$ vectors. 

Crucially, all operations in Algorithm \ref{alg:algorithm} can be written as fully vectorized code, including Line \ref{eq:Qhat}'s maximization along an axis. Only the time stages $t$ require a sequential for-loop. Thus, by design, Algorithm \ref{alg:algorithm} is a natural setting for highly parallelized single instruction, multiple data (SIMD) algebraic operations, for which GPU computation excels.

\section{Data}
\label{sec:data}
Historical day-ahead (DA) and real-time (RT) prices are downloaded from ISO-NE (``Kendall'' price node) for the years 2023-2024. For simplicity we consider hourly averages of RT prices. While our DP method allows any probabilistic forecast of $\{\lambda_t\}_t$, in this work we employ a simple data-driven quantile forecast approach, which uses Y2023 for ``training'' and Y2024 for simulations. For each RT forecast hour in Y2024, we add the corresponding DA price (which is already available) onto 200 samples of empirical quantiles of (RT${-}$DA) spreads from Y2023, conditioned on the same (month, hour), in order to obtain probabilistic forecasts of the RT price.

\section{Numerical results}
We assume 85\% roundtrip efficiency ($\eta = \sqrt{0.85}$), $\overline{p}$ = 1MW power capacity, and $s_0 = 0$MWh initial state, unless otherwise noted. \emph{\textbf{Software}}: LP and MILP models are written in \texttt{CVXPY} and solved with \texttt{Gurobi}; solve times are reported for the time spent inside the solver (which excludes \texttt{CVXPY} compilation time). The discretized DP method is implemented in Python using \texttt{NumPy} for CPU and \texttt{CuPy} \cite{cupy} for GPU. \texttt{CuPy} uses identical software syntax as \texttt{NumPy}, enabling a controlled study of the impact from CPU vs. GPU hardware. \emph{\textbf{Hardware resources}}: Computations are executed on one node of the MIT Engaging high-performance computing cluster. CPU computations are on an Intel Xeon Platinum 8562Y+ CPU processor (32 physical cores, 2.8GHz), GPU computations use an NVIDIA L40S GPU (91.6 TFLOPS for FP32).

\subsection{Accuracy of discretized DP (CPU)}
Imposing a discretization scheme may raise concerns about how accuracy depends on the chosen discretization granularity. So we first benchmark our implemented DP method (still on CPU) relative to an exact LP formulation in the deterministic perfect foresight setting, solving with the full 8784-hour RT prices from Y2024 as $\{\lambda_t\}$ with a single stochastic sample, $R=1$. The problem is solved for a battery with a 4-hour duration, and the DP method is tested at increasingly granular levels of $\delta$. We also obtain bids as in Eq. \eqref{eq:bid_curves}, and evaluate them using a merit-order market clearing simulation based on the realized prices. That is, assuming the storage is a price-taker, for each $t$ intersect the realized RT price with the bid curve to obtain the cleared quantity award $p_t$. 

\begin{table}[ht]
\centering
\caption{Deterministic prices (8784 hours): Accuracy of discretized DP}
\label{tab:opt_comparison}
\setlength{\tabcolsep}{3pt}
\begin{tabular}{c ccccccc}
\toprule
\multirow{2}{*}{\textbf{Method}} & \multirow{2}{*}{$\delta$} & \textbf{Actions} & \textbf{Time} & \multicolumn{2}{c}{\textbf{Quantity opt.}} & \multicolumn{2}{c}{\textbf{Bidding}} \\
 &  & ($P$) & (sec) & Obj. (k\$) & Gap & Obj. (k\$) & Gap \\
\midrule
LP & -- & -- & 0.25 & 57.65 & -- & -- & -- \\
\midrule
\multirow{4}{*}{DP} & 0.10 & 22 & 0.15 & 57.54 & $-0.19\%$ & 57.55 & $-0.17\%$ \\
 & 0.05 & 42 & 0.23 & 57.57 & $-0.13\%$ & 57.59 & $-0.09\%$ \\
 & 0.02 & 103 & 0.64 & 57.63 & $-0.04\%$ & 57.63 & $-0.03\%$ \\
 & 0.01 & 203 & 2.14 & 57.63 & $-0.02\%$ & 57.64 & $-0.02\%$ \\
\bottomrule
\end{tabular}
\end{table}

Table \ref{tab:opt_comparison} reports $\delta$ resolutions and action space sizes $P$, with solve times in seconds. Objectives from quantity optimization (``Quantity opt.'') and price-quantity bidding (``Bidding'') are reported in thousand USD per MW-year, along with relative optimality gaps compared to the exact LP. The DP objectives remain below the LP upper bound, which numerically verifies that Section \ref{sec:discretization}'s discretization framework indeed results in a restriction of the original continuous problem. A 0.2\% accuracy is already achieved with a fairly coarse discretization (22 power actions), suggesting that such a resolution is sufficient for practical valuation and bidding purposes. Accuracy further improves at higher resolutions, to within 0.02\%, suggesting that accuracy scales roughly linearly with $\delta$.

\subsection{Case study of negative prices (CPU)}

A limitation of the prior analytic DP approach in \cite{xu2020operational} is the ex-ante prevention of discharging during negative price periods, as well as the requirement of linear (or at least concave) payoff functions. To demonstrate the generalized discretized DP approach's economic value, we create an artificial price time series by level-shifting the prices to create a deterministic price series consisting of all negative prices:
$$\lambda_t^{neg} = \lambda_t - \left(\max_t \lambda_t\right) \leq 0, \quad \forall t \in [T].$$

While artificially constructed, a setting of frequent negative pricing is generally plausible, e.g. at a highly congested node. Table \ref{tab:negative_comparison} is a case study to optimize the first 72 hours of the adjusted Y2024 prices series $\{\lambda_t^{neg}\}$, with initial $s_0 = \overline{s}$. We compare our discretized DP method (along with its produced bid curves) with a set of different energy arbitrage formulations. \emph{\textbf{MILP exact}}: uses Eq. \eqref{eq:arbitrage} and \eqref{eq:MILP}. \emph{\textbf{LP relaxation}}: replaces Eq. \eqref{eq:MILP:3d} with a continuous $z_t \in [0,1]$. \emph{\textbf{LP restriction}}: prohibits discharge during negative prices, i.e.
$$p^d_t = 0, \quad \forall t \in [T]: \ \lambda_t \le 0,$$
which replicates the same restriction in the analytic DP method of   \cite{xu2020operational}. Ref. \cite{xu2020operational} explains how this restriction is sufficient to eliminate simultaneous charge-discharge. Thus, this is a valid restriction of the MILP exact feasible region. \textbf{DP}: our DP method solved with CPU, and then using the convexification approach to produce monotonic bids as in Eq. \eqref{eq:bids_convex}. 

\begin{table}[ht]
\centering
\caption{Negative prices (72 hours): Accuracy of discretized DP, starting at $\overline{s}$}
\label{tab:negative_comparison}
\setlength{\tabcolsep}{3pt}
\resizebox{\columnwidth}{!}{%
\begin{tabular}{
    c
    S[table-format=2.3]
    S[table-format=2.3]
    S[table-format=3.2]
    S[table-format=2.0]
}
\toprule
\textbf{Method} & {\textbf{Time} (sec)} & {\textbf{Obj.} (k\$)} & \textbf{\ Gap} & {\textbf{Freq.\ of} $p^c_t p^d_t > 0$} \\
\midrule
MILP exact     & 25.068 & 12.811 & {\ --} & 0\%  \\
LP relaxation  & 0.071 & 23.045 & 79.89\% & 86\% \\
LP restriction & 0.142 & 0.000  & -100.00\% & 0\% \\
\midrule
\multirow{2}{*}{DP} & 0.003 & 12.777 & -0.27\% & 0\% \\
      & {(Bidding)}    & 12.799 & -0.10\% & 0\% \\
\bottomrule
\end{tabular}%
}
\end{table}

Table \ref{tab:negative_comparison} demonstrates that the LP relaxation can significantly overstate storage profitability, as a consequence of very frequent complementarity violations (``Frequency of $p_t^c p_t^d >0$''). Conversely, the LP restriction approach earns exactly 0 profit when starting with the full $s_0{\ =\ }\overline{s}$, because the storage device is always prohibited to discharge in this case. In contrast, our discretized DP approach guarantees feasible dispatch, since simultaneous charge-discharge is impossible with only one net $p_t$ variable per $t$. Furthermore, our DP method achieves profits that essentially match the exact MILP solution within $0.1{\sim}0.3\%$ accuracy, while being ${\sim}$8,000x faster than MILP.

\subsection{GPU acceleration with multistage stochastic prices}

Table \ref{tab:gpu} tests the DP method in the multistage stochastic price setting, with 200 price samples per time period. A range of different storage durations (energy-to-power ratio, measured in hours) is tested, at two different $\delta$ resolutions. The table lists the discretization granularity in terms of number of $S$ states and $P$ actions. Times to solve the full-year problem, using CPU (\texttt{NumPy}) versus GPU (\texttt{CuPy}) computation, are reported in seconds. As problem size scales up, either in terms of granularity or storage duration, the GPU times remain the same or modestly increase, which significantly outperforms the CPU's solve times by up to about 100x. 

\begin{table}[ht]
\centering
\caption{CPU vs GPU solve times (8784 hours, 200 price samples per hour)}
\begin{tabular}{
        S[table-format=3.0]
        S[table-format=5.0]
        S[table-format=3.0]
        S[table-format=4.2]
        S[table-format=2.2]
        S[table-format=2.1]}
\toprule
\textbf{Duration} & \textbf{States} & \textbf{Actions} & \textbf{CPU time}  & \textbf{GPU time} & \textbf{\emph{GPU}}
\\
{(hours)} & {($S$)} & {($P$)} & {(sec)} & {(sec)} & {\textbf{\emph{speedup}}}
\\
\midrule
4   & 41    & 22  & 3.98    & 2.77   & 1.4{x} \\
20  & 201   & 22  & 18.29   & 2.81   & 6.5{x} \\
100 & 1001  & 22  & 98.12   & 2.94   & 33.4{x} \\
\midrule
4   & 401   & 203 & 397.80  & 4.88   & 81.5{x} \\
20  & 2001  & 203 & 1996.02 & 19.17  & 104.1{x} \\
100 & 10001 & 203 & 9972.64 & 86.67  & 115.1{x} \\
\bottomrule
\end{tabular}
\label{tab:gpu}
\end{table}

\vspace{-5px}
\subsection{Quantifying the value of stochastic DP bid curves (GPU)}

We study the interaction between price uncertainty and storage duration. We compare a range of strategies representing decreasing levels of information access. \emph{\textbf{Perfect foresight}}: LP on realized RT prices (giving an upper bound on profits). \emph{\textbf{Stochastic DP bid curves}}: DP using Section \ref{sec:data}'s stochastic prices, producing price-quantity bid curves as in Eq. \eqref{eq:bid_curves}. \emph{\textbf{Stochastic DP self-scheduled}}: the same stochastic DP model, but we intersect each bid curve for time $t$ with the realized 1-hour-lagged RT price (as an observable price forecast) in order to obtain a single self-scheduled quantity $p_t$. \emph{\textbf{Myopic}}: deterministic LP on DA prices to obtain self-scheduled quantities, which are multiplied by RT prices. 

\begin{figure}[ht]
    \centering
    \includegraphics[width=0.48\linewidth]{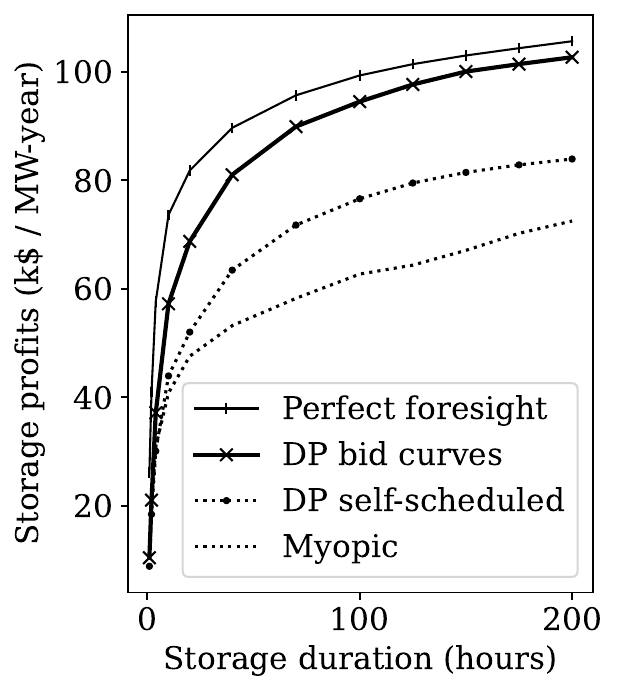}
    \includegraphics[width=0.48\linewidth]{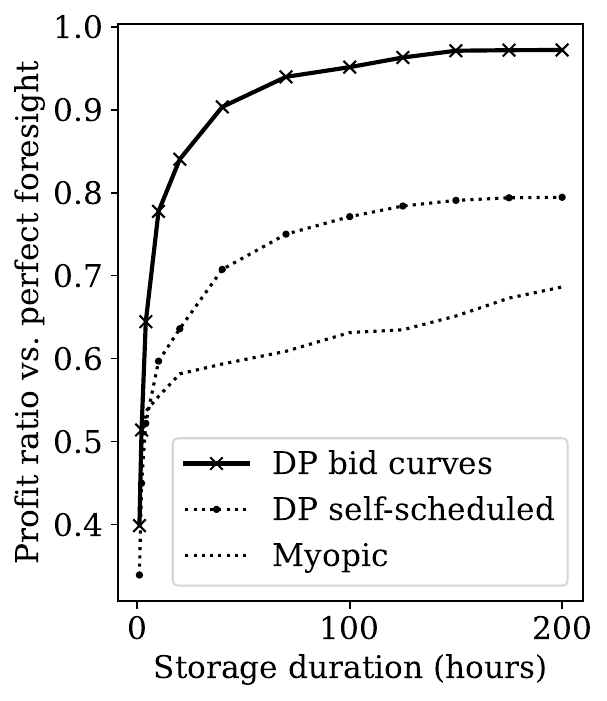}
    \caption{Impact of duration and dispatch strategy on energy arbitrage profits. ``\emph{DP bid curves}'' achieve up to 32\% higher profits than ``\emph{DP self-scheduled}'' (this difference represents the value of using price-quantity bid curves), which in turn achieves up to 24\% higher profits than ``\emph{Myopic}''.}
    \label{fig:profits}
\end{figure}

As seen in Fig. \ref{fig:profits}, across a range of storage durations, using stochastic modeling significantly outperforms the myopic dispatch (by up to 24\%), which represents the value of incorporating uncertainty. Furthermore, utilizing bidding curves rather than single quantities leads to additional gains in profitability (by up to 32\%). An intuition for this additional value is that bid curves can embody the price-contingent optimal policy function. As duration increases, the realizable arbitrage profits also increase as a percentage capture of the perfect foresight upper bound. The evolution of these profit capture ratios has crucial implications for energy storage investment decisions.

Solving all DP models for Fig. \ref{fig:profits} takes 35 seconds with GPU acceleration, compared to about 866 seconds with CPU (14 minutes, or 25x longer). This speedup highlights the power of the proposed method to utilize GPU acceleration in an efficient, on-the-fly manner without having to ``retrain'', unlocking wider studies across practically relevant scenarios and parameters.

\vspace{-5px}
\section{Conclusion}

This paper applies GPU acceleration to dynamic programming in order to speed up calculations for energy arbitrage with multistage price uncertainty. This is made possible by our proposed tensor-based algorithmic formulation of DP backward induction. Up to about 100x speedup is found for GPU vs. CPU computation.

A suite of experiments validate the numerical accuracy of our proposed DP method. Compared to an existing analytic DP approach's convex restriction, our discretized method significantly improves profitability under scenarios with frequent negative prices, while being 8,000x faster than an exact MILP solver. We also demonstrate the ability of the DP method to generate price-quantity bid curves that meaningfully outperform quantity-only dispatch.

The GPU-accelerated DP method can enable more effective investigation over a wide range of energy storage parameters, as well as different price forecasting methods, which can be explored in future work. This approach enables fast and accurate storage valuations, providing a practical tool for large-scale scenario analysis and storage investment studies.

\bibliographystyle{IEEEtran}
\bibliography{ref}

\newpage

\vfill

\end{document}